\documentclass[leqno,12pt]{article}

\usepackage{amsmath,amssymb,amsthm,graphics}

\makeatletter 
\newif\ify@autoscale \y@autoscaletrue \def\Yautoscale#1{\ifnum #1=0 
  \y@autoscalefalse\else\y@autoscaletrue\fi} 
\newdimen\y@b@xdim 
\newdimen\y@boxdim \y@boxdim=13pt 
\def\Yboxdim#1{\y@autoscalefalse\y@boxdim=#1} 
\newdimen\y@linethick    \y@linethick=.3pt 
\def\Ylinethick#1{\y@linethick=#1} 
\newskip\y@interspace \y@interspace=0ex plus 0.3ex 
\def\Yinterspace#1{\y@interspace=#1} 
\newif\ify@vcenter   \y@vcenterfalse 
\def\Yvcentermath#1{\ifnum #1=0 \y@vcenterfalse\else\y@vcentertrue\fi} 
\newif\ify@stdtext   \y@stdtextfalse 
\def\Ystdtext#1{\ifnum #1=0 \y@stdtextfalse\else\y@stdtexttrue\fi} 
\newif\ify@enable@skew   \y@enable@skewfalse 
\def\y@vr{\vrule height0.8\y@b@xdim width\y@linethick depth 0.2\y@b@xdim} 
\def\y@emptybox{\y@vr\hbox to \y@b@xdim{\hfil}} 
\ify@enable@skew 
 \def\y@abcbox#1{\if :#1\else 
   \y@vr\hbox to \y@b@xdim{\hfil#1\hfil}\fi} 
 \def\y@mathabcbox#1{\if :#1\else 
   \y@vr\hbox to \y@b@xdim{\hfil$#1$\hfil}\fi} 
\else 
 \def\y@abcbox#1{\y@vr\hbox to \y@b@xdim{\hfil#1\hfil}} 
 \def\y@mathabcbox#1{\y@vr\hbox to \y@b@xdim{\hfil$#1$\hfil}} 
\fi 
\def\y@setdim{%
  \ify@autoscale%
   \ifvoid1\else\typeout{Package youngtab: box1 not free! Expect an 
     error!}\fi%
   \setbox1=\hbox{A}\y@b@xdim=1.6\ht1 \setbox1=\hbox{}\box1%
  \else\y@b@xdim=\y@boxdim \advance\y@b@xdim by -2\y@linethick 
  \fi} 
\newcount\y@counter 
\newif\ify@islastarg 
\def\y@lastargtest#1,#2 {\if\space #2 \y@islastargtrue 
  \else\y@islastargfalse\fi} 
\def\y@emptyboxes#1{\y@counter=#1\loop\ifnum\y@counter>0 
  \advance\y@counter by -1 \y@emptybox\repeat} 
\def\y@nelineemptyboxes#1{%
  \vbox{%
    \hrule height\y@linethick%
    \hbox{\y@emptyboxes{#1}\y@vr} 
    \hrule height\y@linethick}\vspace{-\y@linethick}} 
\def\yng(#1){%
  \y@setdim%
  \hspace{\y@interspace}%
  \ifmmode\ify@vcenter\vcenter\fi\fi{%
  \y@lastargtest#1, 
  \vbox{\offinterlineskip 
    \ify@islastarg 
     \y@nelineemptyboxes{#1} 
    \else 
     \y@ungempty(#1) 
    \fi}}\hspace{\y@interspace}} 
\def\y@ungempty(#1,#2){%
  \y@nelineemptyboxes{#1} 
  \y@lastargtest#2, 
  \ify@islastarg 
   \y@nelineemptyboxes{#2} 
  \else 
   \y@ungempty(#2) 
  \fi} 
\def\y@nelettertest#1#2. {\if\space #2 \y@islastargtrue 
  \else\y@islastargfalse\fi} 
\def\y@abcboxes#1#2.{%
  \ify@stdtext\y@abcbox#1\else\y@mathabcbox#1\fi%
  \y@nelettertest #2. 
  \ify@islastarg\unskip%
   \ify@stdtext\y@abcbox{#2}\else\y@mathabcbox{#2}\fi%
  \else\y@abcboxes#2.\fi} 
\ify@enable@skew 
 \newdimen\y@full@b@xdim 
 \newcount\y@m@veright@cnt 
 \def\y@get@m@veright@cnt#1#2.{%
   \if :#1 \advance\y@m@veright@cnt by 1\y@get@m@veright@cnt#2.\fi} 
 \let\y@setdim@=\y@setdim 
 \def\y@setdim{%
   \y@setdim@ \y@full@b@xdim=\y@b@xdim 
   \advance\y@full@b@xdim by 1\y@linethick} 
 \def\y@m@veright@ifskew#1{ 
   \y@m@veright@cnt=0 \y@get@m@veright@cnt#1. 
   \moveright \y@m@veright@cnt\y@full@b@xdim} 
\else 
 \def\y@m@veright@ifskew#1{} 
\fi 
\def\y@nelineabcboxes#1{%
  \y@nelettertest #1. 
  \ify@islastarg 
   \y@m@veright@ifskew{#1} 
    \vbox{ 
      \hrule height\y@linethick%
      \hbox{\ify@stdtext\y@abcbox#1\else\y@mathabcbox#1\fi\y@vr} 
      \hrule height\y@linethick}\vspace{-\y@linethick} 
  \else 
   \y@m@veright@ifskew{#1} 
    \vbox{ 
      \hrule height\y@linethick%
      \hbox{\y@abcboxes #1.\y@vr}%
      \hrule height\y@linethick}\vspace{-\y@linethick} 
  \fi} 
\def\young(#1){%
  \y@setdim%
  \hspace{\y@interspace}%
  \y@lastargtest#1, 
  \ifmmode\ify@vcenter\vcenter\fi\fi{%
  \vbox{\offinterlineskip 
    \ify@islastarg\y@nelineabcboxes{#1}%
    \else\y@ungabc(#1)%
    \fi}}\hspace{\y@interspace}} 
\def\y@ungabc(#1,#2){%
  \y@nelineabcboxes{#1}%
  \y@lastargtest#2, 
  \ify@islastarg\y@nelineabcboxes{#2}%
  \else\y@ungabc(#2)%
  \fi} 
\y@vcentertrue 
\makeatother

\newcommand{\g}{\mathfrak{g}}

\newcommand{\Da}{\mathfrak{a}}
\newcommand{\p}{\mathfrak{p}}
\newcommand{\n}{\mathfrak{n}}
\newcommand{\Prim}{\mathrm{Prim}}
\newcommand{\Conj}{\mathrm{Conj}}
\newcommand{\vol}{\mathrm{vol}}

\newcommand{\CSR}{\mathrm{CSR}}
\newcommand{\Id}{\mathrm{Id}}
\newcommand{\Ad}{\mathrm{Ad}}

\newcommand{\tr}{\mathrm{tr}}
\newcommand{\li}{\mathrm{li}}
\newcommand{\Ind}{\mathrm{Ind}}

\newcommand{\as}{\quad\text{as}\quad}
\newcommand{\tinf}{\to\infty}
\newcommand{\disp}{\displaystyle}
\newcommand{\bsla}{\backslash}
\newcommand{\nt}{\notag}

\newcommand{\upar}{\uparrow}

\newcommand{\bC}{\mathbb{C}}
\newcommand{\bR}{\mathbb{R}}
\newcommand{\bQ}{\mathbb{Q}}
\newcommand{\bZ}{\mathbb{Z}}
\newcommand{\bH}{\mathbb{H}}
\newcommand{\bN}{\mathbb{N}}
\newcommand{\noi}{\noindent}
\newcommand{\divset}{\hspace{3pt}|\hspace{3pt}}
\newcommand{\bigdivset}{\hspace{3pt}\big|\hspace{3pt}}
\newcommand{\Bigdivset}{\hspace{3pt}\Big|\hspace{3pt}}

\newcommand{\Biggdivset}{\hspace{3pt}\Bigg|\hspace{3pt}}

\newcommand{\Gam}{\Gamma}

\newcommand{\sr}{\mathrm{SL}_2(\bR)}
\newcommand{\sz}{\mathrm{SL}_2(\bZ)}

\newtheorem{thm}{Theorem}[section]
\newtheorem{prop}[thm]{Proposition}
\newtheorem{lem}[thm]{Lemma}
\newtheorem{cor}[thm]{Corollary}
\newtheorem{rem}[thm]{Remark}
\newtheorem{prob}[thm]{Problem}

\def\miniyng(#1){\scalebox{0.4}{$\yng(#1)$}}

\numberwithin{equation}{section}

\title{Splitting density for lifting about discrete groups}
\author{\textsc{Yasufumi Hashimoto}\thanks{Supported by 
JSPS Research Fellowships for Young Scientists.}
\textsc{ and Masato Wakayama}\thanks{Partially supported by JSPS Grant-in-Aid for Scientific Research 
(B) No. 15340012.}}
\date{}

\begin{document}
\markboth
{Y. Hashimoto and M. Wakayama}
{Splitting density for lifting about discrete groups}
\pagestyle{myheadings}

\maketitle 


\begin{abstract}
We study splitting densities of primitive elements of a discrete subgroup 
of  a connected non-compact semisimple Lie group 
of real rank one with finite center in another larger such discrete subgroup.
When the corresponding cover of such a locally symmetric negatively curved Riemannian manifold 
is regular, 
the densities can be easily obtained from the results due to Sarnak or Sunada. 
Our main interest is a case where the covering is not necessarily regular.
Specifically, for the cases of the modular group and its congruence subgroups,  
we determine the splitting densities explicitly. 
As an application, we study analytic properties of the zeta function defined by 
the Euler product over elements consisting all primitive elements 
which satisfy a certain splitting law for a given lifting.
\end{abstract}

\section{Introduction}
Let $\bH$ be the upper half plane and 
$\Gamma$ a discrete subgroup of $\sr$ such that $\vol(\Gamma\bsla\bH)<\infty$.
Let $\Prim(\Gamma)$ be the set of primitive hyperbolic conjugacy classes of $\Gamma$, 
$N(\gamma)$ the square of the larger eigenvalue of $\gamma\in\Prim(\Gamma)$, 
and $\pi_{\Gamma}(x)$ the number of $\gamma\in\Prim(\Gamma)$ satisfying $N(\gamma)<x$. 
Then, the so-called prime geodesic theorem for $\Gamma$ was discovered by Selberg \cite{Se}
in the early 1950s.
In fact, it tells now (see also \cite{Sa} and \cite{He}) that
\begin{align}
\pi_{\Gamma}(x)=\li{(x)}+O(x^{\delta})\as x\tinf,\label{pgt}
\end{align}
where $\li{(x)}:=\int_{2}^{x}(1/\log{t})dt$ and the constant $\delta$ $(0<\delta<1)$ depends 
on $\Gamma$. 
By using the one-to-one correspondence due to Gauss \cite{G}
between the primitive hyperbolic conjugacy classes of $\sz$ and 
the equivalence classes of the primitive indefinite binary quadratic forms, 
inspired by the work of Selberg, 
Sarnak \cite{Sa} obtained an asymptotic behavior of the sum of the class numbers 
of the quadratic forms from the prime geodesic theorem for $\Gamma=\sz$.
A certain extension of the result in \cite{Sa} for congruence subgroups of $\sz$ 
was recently made in Hashimoto \cite{H}.

The aim of this paper is to study various splitting densities of primitive elements 
of $\tilde{\Gamma}$ in $\Gamma$, where $\tilde{\Gamma}$ denotes 
a subgroup of $\Gamma$ of finite index. 
Although the elements of $\tilde{\Gamma}$ are obviously those of $\Gamma$, 
primitive elements of $\tilde{\Gamma}$ are not necessarily primitive in $\Gamma$.
We consider a problem asking how many primitive elements of $\tilde{\Gamma}$ remain 
also primitive in $\Gamma$,
and moreover, how many primitive elements of $\tilde{\Gamma}$ which are not primitive in $\Gamma$ 
are equal to a given power of primitive elements of $\Gamma$? 

Historically, this kind of branching problem is quite fundamental in 
algebraic number theory. 
Actually, for algebraic extensions of algebraic number fields, 
similar problems had been studied 
by, for example, Artin \cite{Ar}, Tchebotarev \cite{Tc} and Takagi \cite{Ta} 
in the early 20th century. 
The problem for algebraic number fields can be drawn as follows;
let $k$ be an algebraic number field over $\bQ$ and 
$K$ a algebraic extension of $k$ with $n:=[K:k]<\infty$.
We denote by $N_k(\Da)$ the norm of an ideal $\Da$ in $k$.
For a given prime ideal $\p$ of $k$ unramified in $K$, 
there exist a finite number of prime ideals $\p_1,\cdots,\p_m$ of $K$ and
positive integers $e_1,\cdots,e_m$ ($e_1\geq\cdots e_m\geq1$) such that 
$\p=\p_1\cdots\p_m$ and $N_K(\p_i)=N_k(\p)^{e_i}$. 
Since the sum $\sum_{i=1}^{m}e_i$ equals $n$, $(e_1,\cdots,e_m)$ is a partition of $n$.
We call that a prime ideal $\p$ of $k$ is $\lambda$-type in $K$ 
when $\lambda=(e_1,\cdots,e_m)\vdash n$. 
What is the main question is, for a given $\lambda\vdash n$, to count the number 
of prime ideals of $k$ which are $\lambda$-type in $K$.

We now formulate our problem in terms of the geometry 
of negatively curved locally symmetric Riemannian manifolds, 
precisely, by use of the lifting of the primitive geodesics.
First, we prepare some notations.

Let $G$ be a connected non-compact semisimple Lie group of real rank one with finite center
and $G=KA_{\p}N$ be an Iwasawa decomposition of $G$. 
We denote by $\g,\Da_{\p},\n$ the Lie algebra of $G,A_{\p},N$ respectively. 
For the Cartan involution $\theta$ of $G$, $\Da\supset\Da_{\p}$ is defined 
as a $\theta$-stable Cartan subalgebra of $\g$.
Let $\g^{\bC},\Da^{\bC}$ be the complexifications of $\g,\Da$ respectively.  
We denote by $\Phi^{+}$ an $\Da_{\p}$-compatible system of positive roots 
in the set of nonzero roots of $(\g^{\bC},\Da^{\bC})$, 
$P^{+}=\{\alpha\in\Phi^{+}\divset\text{$\alpha\not\equiv0$ on $\Da_{\p}$}\}$ and
$\Sigma^{+}$ the set of the restrictions of the elements of $P^{+}$ on $\Da_{\p}$. 
Then $\Sigma^{+}$ is written as $\Sigma^{+}=\{\beta\}$ or $\{\beta,2\beta\}$ for some $\beta$. 
We choose $H_0\in \Da_{\p}$ such that $\beta(H_0)=1$ and put 
$\rho=1/2\sum_{\alpha\in P^{+}}\alpha$ and $\rho_0=\rho(H_0)$. 

Let $\Gamma$ be a discrete subgroup of $G$ such that the volume 
of $X_{\Gamma}:=\Gamma\bsla G/K$ is finite. 
We denote by $\Prim(\Gamma)$ a set of primitive hyperbolic conjugacy classes of $\Gamma$. 
For $\gamma\in\Gamma$, the norm $N(\gamma)$ is defined by  
\begin{align*}
N(\gamma)=\max\{|\delta|^k\divset\text{$\delta$ is an eigenvalue of $\Ad(\gamma)$}\},
\end{align*}
where $\Ad$ is the adjoint representation of $G_{\bC}$, the analytic group 
with Lie algebra $\g_{\bC}$ 
and $k(=1,2)$ denotes the number of elements in $\Sigma^{+}$.
Denote $\pi_{\Gamma}(x)$ by the number of $\gamma\in \Prim(\Gamma)$ satisfying $N(\gamma)<x$. 
Then $\pi_{\Gamma}(x)$ behaves
\begin{align*}
\pi_{\Gamma}(x)=\li(x^{2\rho_0})+O(x^{\delta})\as x \tinf,
\end{align*}
where $\delta(0<\delta<2\rho_0)$ is a constant depending on $\Gamma$ (see, e.g. \cite{GW}). 

\begin{prob}\label{prob}
Let $\tilde{\Gamma}$ be a subgroup of $\Gamma$ of finite index and 
suppose that $X_{\tilde{\Gamma}}$ is a finite cover of $X_{\Gamma}$.
We denote by $p$ a natural projection from $X_{\tilde{\Gamma}}$ to $X_{\Gamma}$. 
Let $C_{\gamma}$ be a closed primitive geodesic of $X_{\Gamma}$ corresponding 
to $\gamma\in\Prim(\Gamma)$, 
and $l(\gamma)$ the length of $C_{\gamma}$ ($N(\gamma):=e^{l(\gamma)}$). 
For a given $\gamma\in\Prim(\Gamma)$, 
there exists a finite number of elements $\gamma_1,\gamma_2,\cdots,\gamma_k$ 
of $\Prim(\tilde{\Gamma})$ 
and positive integers $m_1,\cdots,m_k$ such that $p(C_{\gamma_j})=C_{\gamma}$ 
with $l(\gamma_j)=m_j l(\gamma)$. 
We may assume that $m_1\geq\cdots\geq m_k$.
Since $\sum_{j=1}^{k}m_j=[\Gamma:\tilde{\Gamma}](=:n)$, 
$(m_1,m_2,\cdots,m_k)$ is considered as a partition of $n$. 
We call an element $\gamma\in\Prim(\Gamma)$ is $\lambda$-type in $\tilde{\Gamma}$ 
when $\lambda=(m_1,m_2,\cdots,m_k)\vdash n$.
We define $\pi_{\tilde{\Gamma}\upar\Gamma}^{\lambda}(x)$ and 
its density $\mu_{\tilde{\Gamma}\upar\Gamma}^{\lambda}(x)$ relative to $\pi_{\Gamma}(x)$ as 
\begin{align*}
\pi_{\tilde{\Gamma}\upar\Gamma}^{\lambda}(x):=&\#\{\gamma\in\Prim(\Gamma)\divset
\text{$\gamma$ is $\lambda$-type in $\tilde{\Gamma}$},N(\gamma)<x\},\\
\mu_{\tilde{\Gamma}\upar\Gamma}^{\lambda}(x)
:=&\pi_{\tilde{\Gamma}\upar\Gamma}^{\lambda}(x)/\pi_{\Gamma}(x).
\end{align*}
For a given $\lambda\vdash n$, study the asymptotic behavior of 
$\pi_{\tilde{\Gamma}\upar\Gamma}^{\lambda}(x)$ and 
$\mu_{\tilde{\Gamma}\upar\Gamma}^{\lambda}(x)$ when $x\tinf$.\qed
\end{prob}

If a covering $X_{\tilde{\Gamma}}\to X_{\Gamma}$ is regular, that is, 
$\tilde{\Gamma}$ is a normal subgroup of $\Gamma$, 
the problem can be easily solved (Theorem \ref{thm2}) based on the results 
in \cite{Sa} and \cite{Su2}.
Thus, the main focus of the present paper is a study of Problem \ref{prob} 
when $X_{\tilde{\Gamma}}$ is not necessarily a regular cover of $X_{\Gamma}$. 
Especially, in the case where $G=\sr$, $\Gamma=\sz$ and $\tilde{\Gamma}$ 
is a congruence subgroup of $\sz$, 
the splitting densities can be obtained explicitly (see Section 4 and 5).

Applying the results in Section 4 and 5 to Venkov-Zograf's formula \cite{VZ} about 
the relation between the Selberg zeta functions for $\Gam$ and $\tilde{\Gam}$, 
we can obtain an expression of the Selberg zeta function for the congruence subgroup 
as a product over elements of $\Prim(\sz)$. 
Then, in the last section, by taking a quotient of such expressions of Selberg's zeta functions 
for two congruence subgroups, 
we give a functional equation and an analytic continuation to the right half plane 
of the zeta function defined by the Euler product over elements 
consisting all primitive elements of $\Prim(\sz)$
which satisfy a certain splitting law for a given lifting in the congruence subgroup.

\section{General cases}
It is not true in general that $\pi_{\tilde{\Gamma}\upar\Gamma}^{\lambda}(x)>0$ for 
$\lambda\vdash n=[\Gamma:\tilde{\Gamma}]$ as we see below. 
Actually, $\pi_{\tilde{\Gamma}\upar\Gamma}^{\lambda}(x)=0$ may hold 
for many partitions $\lambda\vdash n$.
Hence, in Problem \ref{prob}, it is important to determine 
partitions $\lambda$ of which the density $\pi_{\tilde{\Gamma}\upar\Gamma}^{\lambda}(x)$ 
is positive.
For a general pair $(\Gamma,\tilde{\Gamma})$ such that $\tilde{\Gamma}\subset \Gamma$, 
we have the following basic theorem.

\begin{thm}\label{thm1}
Let $\Gamma'$ be the (unique) maximal normal subgroup of $\Gamma$ contained in $\tilde{\Gamma}$. 
Let $\Xi:=\Gamma/\Gamma'$ and $\Conj(\Xi)$ the set of conjugacy classes of $\Xi$.
We denote by $M(\gamma):=\min\{m\geq1\divset\gamma^m\in\ \Gam'\}$ for $\gamma\in\Gamma$ and
$A_{\Gamma'\upar\Gamma}:=\{M(\gamma)\divset\gamma\in \Gamma\}\subset\bN$. 
Define 
\begin{align*}
\Lambda
:=\{(m_1,m_2,\cdots ,m_k)\vdash n\divset\exists M\in A_{\Gamma'\upar\Gamma}, \forall m_i|M\}.
\end{align*}
Then, for $\lambda\in\Lambda$, we have 
\begin{align*}
\pi_{\tilde{\Gamma}\upar\Gamma}^{\lambda}(x)=
\bigg(\sum_{\begin{subarray}{c}[\gamma]\in\Conj(\Xi),
\\ \text{$[\gamma]$ is $\lambda$-type in $\tilde{\Gamma}$}\end{subarray}}
\frac{\#[\gamma]}{|\Xi|}\bigg)\li{(x^{2\rho_0})}+O(x^{\delta})\as x\tinf.
\end{align*}
For $\lambda\not\in\Lambda$, we have $\pi_{\tilde{\Gamma}\upar\Gamma}^{\lambda}(x)=0$.
\end{thm}

\begin{cor}\label{cor}
We have
\begin{align*}
\lim_{x\tinf}\mu_{\tilde{\Gamma}\upar\Gamma}^{\lambda}(x)=\begin{cases}
\disp\sum_{\begin{subarray}{c}[\gamma]\in\Conj(\Xi),\\ 
\text{$[\gamma]$ is $\lambda$-type in $\tilde{\Gamma}$}\end{subarray}}
\frac{\#[\gamma]}{|\Xi|}&\quad\text{for $\lambda\in\Lambda$},\\
0&\quad\text{for $\lambda\not\in\Lambda$}.\end{cases}
\end{align*}
\end{cor}

To prove Theorem \ref{thm1}, we need some preparations.
\begin{lem}\label{prop}
Let $\Psi:=\tilde{\Gamma}/\Gamma'$ and $\lambda:=(m_1,\cdots,m_k)\vdash n$. 
Denote by $\sigma$ the permutation representation of $\Xi$ on $\Xi/\Psi$ 
($\sigma\cong\Ind_{\tilde{\Gamma}}^{\Gamma}1$).
For $\gamma\in\Prim(\Gamma)$, we call that $\sigma(\gamma)$ is $\lambda$-type 
when $\sigma(\gamma)$ is expressed as 
\begin{align}\label{ind}
\sigma(\gamma)\sim\left(\begin{array}{ccc}
S_{m_1}&\cdots&0\\
\vdots&\ddots&\vdots\\
0&\cdots&S_{m_k}
\end{array}\right), 
\end{align}
where $S_{m_i}$ are $m_i\times m_i$-matrices given by
\begin{align*}
S_{m_i}=\begin{cases}1&(m_i=1),\\
\left(\begin{array}{cccccc}
0&1&\cdots&0\\
\vdots&\vdots&\ddots&\vdots\\
0&0&\cdots&1\\
1&0&\cdots&0
\end{array}\right)&(m_i\geq2).\end{cases}
\end{align*}
Then the following two conditions are equivalent.\\
(i) $\sigma(\gamma)$ is $\lambda$-type.\\
(ii) $\gamma\in\Prim(\Gamma)$ is $\lambda$-type in $\tilde{\Gamma}$.
\end{lem}

\begin{proof}
Let $\CSR[\Gamma/\tilde{\Gamma}]$ be 
a complete system of representatives of $\Gamma/\tilde{\Gamma}$.\\
(I) If (i) holds, then there exist $A_1,\cdots,A_k\subset \CSR[\Gamma/\tilde{\Gamma}]$ 
such that $\coprod_{i=1}^{k}A_i=\CSR[\Gamma/\tilde{\Gamma}]$ and 
$$(a_1^{(i)})^{-1}\gamma a_2^{(i)},\cdots,(a_{m_i-1}^{(i)})^{-1}\gamma a_{m_i}^{(i)}, 
(a_{m_i}^{(i)})^{-1}\gamma a_1^{(i)}\in\tilde{\Gamma}$$
for $A_i=(a_1^{(i)},\cdots,a_{m_i}^{(i)})$. 
Hence, we have
\begin{align*}
(a_1^{(1)})^{-1}\gamma^{m_1} a_1^{(1)},\cdots,
(a_1^{(k)})^{-1}\gamma^{m_k} a_1^{(k)}\in\tilde{\Gamma}.
\end{align*}
If we put $\gamma_i:=(a_1^{(i)})^{-1}\gamma^{m_i} a_1^{(i)}$, it is easy to see that 
$\gamma_i$ is primitive in $\tilde{\Gamma}$ and is $\Gamma$-conjugate to $\gamma^{m_i}$.
Hence, we have $p(C_{\gamma_j})=C_{\gamma}$ and $l(\gamma_j)=m_jl(\gamma)$.\\
(II) If (ii) holds, there exist $b_1,\cdots,b_k\in \CSR[\Gamma/\tilde{\Gamma}]$ 
such that $\gamma_i=b_i^{-1}\gamma^{m_i}b_i$. 
Also, there exists $c_1^{(i)}\in \CSR[\Gamma/\tilde{\Gamma}]$ 
such that $(c_1^{(i)})^{-1}\gamma b_i\in\tilde{\Gamma}$.
Since $b_i^{-1}\gamma^{m_i-1}c_1^{(i)}\in\tilde{\Gamma}$, 
it is easy to see that there exist 
$c_2^{(i)},\cdots,c_{m_i-1}^{(i)}\in \CSR[\Gamma/\tilde{\Gamma}]$ such that 
\begin{align*}
(c_1^{(i)})^{-1}\gamma b_i,(c_2^{(i)})^{-1}\gamma c_1^{(i)},\cdots,b_i^{-1}\gamma c_{m_i-1}^{(i)}
\in\tilde{\Gamma}
\end{align*}
recursively. 
If we assume that there exist $1\leq j_1<j_2\leq m_i-1$ such that $c_{j_1}^{(i)}=c_{j_2}^{(i)}$, 
then $b_i^{-1}\gamma^{m_i-j_2}c_j^{(i)}$,$(c_{j_1}^{(i)})^{-1}\gamma^{j_2-j_1}c_{j_1}^{(i)}$, 
$(c_{j_1}^{(i)})^{-1}\gamma^{j_1}b_i\in\tilde{\Gamma}$. 
It follows hence that $b_i^{-1}\gamma^{m_i-j_2+j_1}b_i\in\tilde{\Gamma}$. 
This contradicts, however, the fact that $\gamma_i=b_i^{-1}\gamma^{m_i}b_i$ 
is primitive in $\tilde{\Gamma}$.
Hence $b_i$'s and $c_j^{(i)}$'s are mutually distinct and 
$\{b_i,c_j^{(i)}\}=\CSR[\Gamma/\tilde{\Gamma}]$. 
Therefore we see that $\sigma(\gamma)$ is $(m_1,\cdots,m_k)$-type.
\end{proof}

By using the trace formula, Sarnak has shown the following analytic distribution.
\begin{prop}\label{sarnak}(Theorem 2.4 in \cite{Sa}, \cite{Su1} or \cite{Su2})
For $[g]\in\Conj(\Xi)$, we have
\begin{align*}
\#\{\gamma\in\Prim(\Gamma)\divset\gamma\Gamma'=[g],N(\gamma)<x\}
\sim\frac{\#[g]}{|\Xi|}\li{(x^{2\rho_0})}+O(x^{\delta})\as x\tinf.\qed
\end{align*}
\end{prop}

\noindent{\bf Proof of Theorem \ref{thm1}.} 
Assume that $\gamma\in\Prim(\Gamma)$ is $(1^{l_1}2^{l_2}\cdots n^{l_n})$-type in $\tilde{\Gamma}$. 
Then, because of Lemma \ref{prop}, $\sigma(\gamma)$ is $(1^{l_1}2^{l_2}\cdots n^{l_n})$-type.
Since $\sigma(\gamma^{M(\gamma)})=\Id$, we have $l_j=0$ for $j\nmid M(\gamma)$.
Hence, for $\lambda\not\in\Lambda$, we have $\mu_{\tilde{\Gamma}\upar\Gamma}^{\lambda}(x)=0$.
Since the type of $\sigma(\gamma)$ is invariant under the $\Xi$-conjugation, 
by Proposition \ref{sarnak}, the desired asymptotic follows. \qed

\section{Regular cover cases}
In Theorem \ref{thm1}, 
we see that $\mu_{\tilde{\Gamma}\upar\Gamma}^{\lambda}(x)=0$ holds for all $\lambda\not\in\Lambda$.
However, even if $\lambda\in\Lambda$, it is not necessarily that 
$\mu_{\tilde{\Gamma}\upar\Gamma}^{\lambda}(x)>0$.
Actually, when $\tilde{\Gamma}$ is a normal subgroup of $\Gamma$ $(\tilde{\Gamma}=\Gamma')$, 
we prove that only rectangle shape partitions can appear non-trivially.

\begin{thm}\label{thm2}
If $\tilde{\Gamma}$ is a normal subgroup of $\Gamma$ ($\tilde{\Gamma}=\Gamma'$),
then, for $\lambda=\lambda(m)=(m^{n/m})$ $(m\in A_{\Gamma'\upar\Gamma})$, we have
\begin{align*}
\pi_{\tilde{\Gamma}\upar\Gamma}^{\lambda}(x)
=\bigg(\sum_{\begin{subarray}{c}[\gamma]\in\Conj(\Xi),\\ 
M([\gamma])=m \end{subarray}}\frac{\# [\gamma]}{|\Xi|}\bigg)\li{(x^{2\rho_0})}+O(x^{\delta})
\as x\tinf.
\end{align*} 
For $\lambda\vdash n$, other than the shape above, 
we have $\pi_{\tilde{\Gamma}\upar\Gamma}^{\lambda}(x)=0$.
\end{thm}

\begin{proof}
Since $\tilde{\Gamma}=\Gamma'$, we have $\Xi/\Psi=\Xi/\{\Id\}=\Xi$. 
For $\gamma\in \Xi$, suppose that there exist $g\in \Xi$ and $l<M(\gamma)$ such that $\gamma^l g=g$.
Then we have $\gamma^l=\Id$, but this contradicts the minimality of $M(\gamma)$. 
Hence, for any $g$, we see that the type of $\sigma(g)$ 
is given as $(m^{n/m})$ $(m\in A_{\Gamma'\upar\Gamma})$.
Consequently, applying Proposition \ref{sarnak}, we obtain the desired result.
\end{proof}

\begin{cor}
We have
\begin{align*}
\lim_{x\tinf}\mu_{\tilde{\Gamma}\upar\Gamma}^{\lambda}(x)
=\begin{cases}
\disp\sum_{\begin{subarray}{c}[\gamma]\in\Conj(\Xi),\\ 
M([\gamma])=m \end{subarray}}\frac{\# [\gamma]}{|\Xi|}
&\quad\text{for $\lambda=\lambda(m)$ $(m\in A_{\Gamma'\upar\Gamma})$},\\
0&\quad\text{otherwise}.\end{cases}
\end{align*}
\end{cor}

\begin{rem}
\rm{In the problem for an algebraic number field, 
similar results to Theorem \ref{thm2} had been obtained.
In fact, for the cases of unramified Galois extensions, 
the corresponding densities are non-zero only when the factorization is of rectangle type 
(see, e.g. \cite{Ar}, \cite{Tc}, \cite{Ta} and \cite{Na}).}
\end{rem}

\section{Cases of congruence subgroups of $\sz$}
Let $N$ be a positive integer. 
In this section, we consider the cases of $\Gamma=\sz$ and $\tilde{\Gamma}$ is 
one of the following congruence subgroups. 
\begin{align*}
\Gamma_0(N):=&\{\gamma\in \sz\divset\gamma_{21}\equiv0\bmod{N}\},\\	
\Gamma_1(N):=&\{\gamma\in \sz\divset\gamma_{11}\equiv\gamma_{22}\equiv\pm1,
\gamma_{21}\equiv0\bmod{N}\},\\	
\Gamma(N):=&\{\gamma\in \sz\divset\gamma_{11}\equiv\gamma_{22}\equiv\pm1,
\gamma_{12}\equiv\gamma_{21}\equiv0\bmod{N}\}.	
\end{align*}

Let $p$ be a prime number and, for simplicity, assume that $p\geq3$.
First, we study the cases of $N=p^r$. 
Note that the maximal normal subgroup of $\Gamma$ contained in $\tilde{\Gamma}$ 
is $\Gamma'=\Gamma(p^r)$, whence $\Xi=\mathrm{SL}_2(\bZ/p^r\bZ)/\{\pm\Id\}$ and
\begin{align*}
|\Xi|=&\frac{1}{2}p^{3r-2}(p^2-1),\\
n=&\begin{cases}p^{r-1}(p+1)&(\tilde{\Gamma}=\Gamma_0(p^r)),\\
\disp\frac{1}{2}p^{2r-2}(p^2-1)&(\tilde{\Gamma}=\Gamma_1(p^r)),\\
\disp\frac{1}{2}p^{3r-2}(p^2-1)&(\tilde{\Gamma}=\Gamma(p^r)).\end{cases}
\end{align*}
In these cases we have the following results.
\begin{thm}\label{congruence}
Let
\begin{align*}
&\lambda_0^{p^r}(1):=(1^{p^{r-1}(p+1)}),\\
&\lambda_0^{p^r}(lp^{k}):=\begin{cases}\Big((l)^{(p^r+p^{r-1}-2)/l},(1)^{2}\Big)& 
(k=0,l|(p-1)/2,l>1),\\
\Big((lp^{k})^{p^{r-k-1}(p-1)/l},(lp^{k-1})^{2p^{r-k-1}(p-1)/l},\cdots \\
 \cdots,(l)^{2(p^{r-k}-1)/l},(1)^{2}\Big)& (k>0,l|(p-1)/2,l>1),\\
\Big((lp^{k})^{p^{r-k-1}(p+1)/l}\Big)& (l|(p+1)/2,l>1),\end{cases}\\
&\lambda_0^{p^r}(p^{k},A):=\Big((p^{k})^{p^{r-k-1}(p-1)},(p^{k-1})^{2p^{r-k-1}(p-1)},
\cdots,(p)^{2p^{r-k-1}(p-1)},(1)^{2p^{r-k}}\Big),\\
&\lambda_0^{p^r}(p^{k},B^{(k)})\\&:=\begin{cases}\Big((p^{k})^{p^{r-k}},(p^{k-2})^{p^{r-k}(p-1)}, 
\cdots, (p^2)^{p^{r-k/2-2}(p-1)},(1)^{p^{r-k/2}}\Big)& (\text{$k$ is even}),\\
\Big((p^{k})^{p^{r-k}},(p^{k-2})^{p^{r-k}(p-1)}, 
\cdots, (p)^{p^{r-(k+3)/2}(p-1)},(1)^{p^{r-(k+1)/2}}\Big)&(\text{$k$ is odd}),
\end{cases}\\
&\lambda_0^{p^r}(p^{k},B^{(m)}):=\Big((p^{k})^{p^{r-k}},(p^{k-2})^{p^{r-k}(p-1)},\cdots,\\
&\cdots, (p^{k-m+2})^{p^{r-k+(m-3)(p-1)/2}},(p^{k-m})^{p^{r-k+(m-1)/2}}\Big)
\quad(1\leq m<k,\text{odd}),\\
&\lambda_0^{p^r}(p^{k},B^{(m,+)}):=
\Big((p^k)^{p^{r-k}},(p^{k-2})^{p^{r-k}(p-1)},\cdots,(p^{k-m+2})^{p^{r-k+m/2-2}(p-1)},\\
&(p^{k-m})^{p^{r-k+m/2-1}(p-2)},(p^{k-m-1})^{2p^{r-k+m/2-1}(p-1)},\cdots,
(p)^{2p^{r-k+m/2-1}(p-1)},(1)^{2p^{r-k+m/2}}\Big)\\
&\quad\quad\quad (\text{$1<m<k$ is even}),\\
&\lambda_0^{p^r}(p^{k},B^{(m,-)}):=\Big((p^{k})^{p^{r-k}},(p^{k-2})^{p^{r-k}(p-1)},\cdots \\
&\cdots, (p^{k-m+3})^{p^{r-k+m/2-3}(p-1)},(p^{k-m+1})^{p^{r-k+m/2-1}}\Big)
\quad (1<m<k,\text{even}),\\
&\lambda_0^{p^r}(p^{k},C):=\Big((p^{k})^{p^{r-k-1}(p+1)}\Big),\\
&\lambda_1^{p^r}(lp^{k}):=\Big((lp^{k})^{p^{2r-k-2}(p^2-1)/2l}\Big)
\quad(l|(p\pm1)/2,l\geq1,0\leq k\leq r-1),\\
&\lambda_1^{p^r}(p^{k},B^{(m)}):=
\Big((p^k)^{p^{2r-k-1}(p-1)/2},(p^{k-1})^{p^{2r-k-2}(p-1)^2/2},\cdots,\\&
\cdots,(p^{k-m+1})^{p^{2r-k-2}(p-1)^2/2},(p^{k-m})^{p^{2r-k-1}(p-1)/2}\Big),\\
&\lambda^{p^r}(l):=\Big(l^{|\Xi|/l}\Big).
\end{align*}
Then we have 
\begin{align*}
&\lim_{x\tinf}\mu_{\Gamma_0(p^r)\upar \sz}^{\lambda_0^{p^r}(1)}(x)
=\lim_{x\tinf}\mu_{\Gamma_1(p^r)\upar \sz}^{\lambda_1^{p^r}(1)}(x)
=\lim_{x\tinf}\mu_{\Gamma(p^r)\upar \sz}^{\lambda^{p^r}(1)}(x)
=\frac{2}{p^{3r-2}(p^2-1)},\\
\end{align*}
\begin{align*}
&\lim_{x\tinf}\mu_{\Gamma_0(p^r)\upar \sz}^{\lambda_0^{p^r}(lp^{k})}(x)
=\lim_{x\tinf}\mu_{\Gamma_1(p^r)\upar \sz}^{\lambda_1^{p^r}(lp^{k})}(x)
=\lim_{x\tinf}\mu_{\Gamma(p^r)\upar \sz}^{\lambda^{p^r}(lp^{k})}(x)\\
&=\begin{cases}\disp\frac{\varphi(l)}{p^{r-1}(p-1)}& (k=0,l|(p-1)/2,l>1),\\
\disp\frac{\varphi(l)}{p^{r-1}(p+1)}& (k=0,l|(p+1)/2,l>1),\\
\disp\frac{\varphi(l)}{p^{r-k}}& (k>0,l|(p-1)/2,l>1),\\
\disp\frac{\varphi(l)(p-1)}{p^{r-k}(p+1)}& (k>0,l|(p+1)/2,l>1),\end{cases}\\
&\lim_{x\tinf}\mu_{\Gamma_0(p^r)\upar \sz}^{\lambda_0^{p^r}(p^{k},A)}(x)
=\frac{1}{p^{3(r-k)}},\\
&\lim_{x\tinf}\mu_{\Gamma_0(p^r)\upar \sz}^{\lambda_0^{p^r}(p^{k},C)}(x)
=\frac{p-1}{p^{3(r-k)}(p+1)},\\
&\lim_{x\tinf}\mu_{\Gamma_1(p^r)\upar \sz}^{\lambda_1^{p^r}(p^{k})}(x)
=\frac{2}{p^{3r-3k-1}(p+1)},\\
&\lim_{x\tinf}\mu_{\Gamma_0(p^r)\upar \sz}^{\lambda_0^{p^r}(p^{k},B^{(m,\pm)})}(x)
=\frac{(p-1)}{p^{3r-3k+m+1}},\\
&\lim_{x\tinf}\mu_{\Gamma_0(p^r)\upar \sz}^{\lambda_0^{p^r}(p^{k},B^{(m)})}(x)
=\lim_{x\tinf}\mu_{\Gamma_1(p^r)\upar \sz}^{\lambda_1^{p^r}(p^{k},B^{(m)})}(x)
=\begin{cases}\disp\frac{2}{p^{3r-2k}}&(m=k),\\ \disp\frac{2(p-1)}{p^{3r-3k+m+1}}& (m<k),
\end{cases}\\
&\lim_{x\tinf}\mu_{\Gamma(p^r)\upar \sz}^{\lambda^{p^r}(p^{k})}(x)
=\begin{cases}\disp\frac{2}{p}&(k=r),\\ \disp\frac{2(p^2+p+1)}{p^{3r-3k+1}(p+1)}&(k<r).\end{cases}
\end{align*}
For any other $\lambda\vdash n$, we have $\pi_{\tilde{\Gamma}\upar\Gamma}^{\lambda}(x)=0$.
\end{thm}

To prove Theorem \ref{congruence}, first, we use the following classification of 
the conjugacy classes of $\mathrm{SL}_2(\bZ/p^r\bZ)/\{\pm I\}$ 
(see \cite{Di} for $r=1$ and \cite{Kl} for larger $r$).
\begin{lem}\label{conjugacy}
Each element of $\mathrm{SL}_2(\bZ/p^r\bZ)/\{\pm I\}-\{I\}$ 
is conjugate to one of the following elements.\\
\begin{align*}
\begin{array}{|c|l|c|c|}
\hline
\gamma&\text{parameter} &M(\gamma)&\#[\gamma]\\
\hline
\begin{pmatrix}\delta &0\\ 0& \delta^{-1}\end{pmatrix}^{sp^{k-1}(p-1)/2l}
&\begin{cases}1\leq k\leq r,\\l|(p-1)/2(l>1),\\s\in(\bZ/lp^{r-k}\bZ)^{*}/\{\pm1\} \end{cases}
&lp^{r-k}&p^{2r-1}(p+1)\\
\hline
\begin{pmatrix}\delta &0\\ 0& \delta^{-1}\end{pmatrix}^{sp^{k-1}(p-1)}
&\begin{cases}1\leq k\leq r-1,\\s\in(\bZ/p^{r-k}\bZ)^{*}/\{\pm1\}\end{cases}&p^{r-k}
&p^{2r-2k-1}(p+1)\\
\hline
\end{array}
\end{align*}
\begin{align*}
\begin{array}{|c|l|c|c|}
\hline
\begin{pmatrix}1+\alpha p^{2k+m} &p^k\\ \alpha p^{k+m}& 1\end{pmatrix}
&\begin{cases}0\leq k\leq r-1,\\ 1\leq m\leq r-k,\\ 
\alpha\in(\bZ/p^{r-k-m}\bZ)^{*}\end{cases}
&p^{r-k}&p^{2r-2k-2}(p^2-1)/2\\
\hline
\begin{pmatrix}1+\nu\alpha p^{2k+m} &\nu p^k\\ \alpha p^{k+m}& 1\end{pmatrix}
&\begin{cases}0\leq k\leq r-1,\\ 1\leq m\leq r-k,\\\alpha\in(\bZ/p^{r-k-m}\bZ)^{*}\end{cases} 
&p^{r-k}&p^{2r-2k-2}(p^2-1)/2\\
\hline
\Omega^{sp^{k-1}(p+1)/2l}
&\begin{cases}1\leq k\leq r,\\l|(p+1)/2,l>1,\\s\in(\bZ/lp^{r-k}\bZ)^{*}/\{\pm1\}\end{cases}
&lp^{r-k}&p^{2r-1}(p-1)\\
\hline
\Omega^{sp^{k-1}(p+1)}
&\begin{cases}1\leq k\leq r-1,\\s\in(\bZ/p^{r-k}\bZ)^{*}/\{\pm1\}\end{cases}
&p^{r-k}&p^{2r-2k-1}(p-1),\\
\hline
\end{array}
\end{align*}
where $\delta$ is a generator of $(\bZ/p^r\bZ)^{*}/\{\pm1\}$, 
$\nu$ is a non-quadratic residue of $p$ 
and $\Omega\in\mathrm{SL}_2(\bZ/p^r\bZ)/\{\pm I\}$ is of order $p^{r-1}(p+1)/2$. 
Note that $(\tr{\Omega})^2-4$ is a non-quadratic residue of $p$. \qed
\end{lem}

The claims for $\tilde{\Gamma}=\Gamma(p^r)$ in Theorem \ref{congruence} follow easily 
from  Theorem \ref{thm2} and Lemma \ref{conjugacy}. 

Now, we put 
\begin{align}\label{conj}
A_0^{(k,l)}:=&\bigcup_{s\in(\bZ/lp^{r-k}\bZ)^{*}/\{\pm1\}}
\bigg[\begin{pmatrix}\delta &0\\ 0& \delta^{-1}\end{pmatrix}^{sp^{k-1}(p-1)/l}\bigg]
\qquad\big(l\big|(p-1)/2\big),\nt\\
A_k:=&\bigcup_{s\in(\bZ/p^{r-k}\bZ)^{*}/\{\pm1\}}
\bigg[\begin{pmatrix}\delta &0\\ 0& \delta^{-1}\end{pmatrix}^{sp^{k-1}(p-1)}\bigg],\nt\\
B_k^{(m)}:=&\bigcup_{\alpha\in(\bZ/p^{r-k-m}\bZ)^{*}}\Bigg(
\bigg[\begin{pmatrix}1+\alpha p^{2k+m} &p^k\\ \alpha p^{k+m}& 1\end{pmatrix}\bigg]\cup
\bigg[\begin{pmatrix}1+\nu\alpha p^{2k+m} &\nu p^k\\ \alpha p^{k+m}& 1\end{pmatrix}\bigg]\Bigg),\\
C_0^{(k,l)}:=&\bigcup_{s\in(\bZ/lp^{r-k}\bZ)^{*}/\{\pm1\}}
\bigg[\Omega^{sp^{k-1}(p+1)/l}\bigg]
\qquad\big(l\big|(p+1)/2\big),\nt\\
C_k:=&\bigcup_{s\in(\bZ/p^{r-k}\bZ)^{*}/\{\pm1\}}
\bigg[\Omega^{sp^{k-1}(p-1)}\bigg].\nt
\end{align}
We divide $B_k^{(m)}$ for even $m<r-k$ by $B_k^{(m)}=B_k^{(m,+)}\cup B_k^{(m,-)}$ where 
\begin{align*}
B_k^{(m,\pm)}:=&\bigcup_{\begin{subarray}{c}\alpha\in(\bZ/p^{r-k-m}\bZ)^{*}\\
(\frac{\alpha}{p})=\pm1\end{subarray}}
\bigg[\begin{pmatrix}1+\alpha p^{2k+m} &p^k\\ \alpha p^{k+m}& 1\end{pmatrix}\bigg]\cup
\bigcup_{\begin{subarray}{c}\alpha\in(\bZ/p^{r-k-m}\bZ)^{*}\\(\frac{\alpha}{p})=\mp1\end{subarray}}
\bigg[\begin{pmatrix}1+\nu\alpha p^{2k+m} &\nu p^k\\ \alpha p^{k+m}& 1\end{pmatrix}\bigg].
\end{align*}
Note that
\begin{align*}
\sz/\Gam(p^r)-\Gamma(p)/\Gam(p^r)
&=\Big(\bigcup_{l|(p-1)/2,l>1}\bigcup_{k=1}^{r}A_0^{(k,l)}\Big)\cup
\Big(\bigcup_{m=1}^{r}B_0^{(m)}\Big)\cup 
\Big(\bigcup_{l|(p+1)/2,l>1}\bigcup_{k=1}^{r}C_0^{(k,l)}\Big),\\
\Gamma(p^k)/\Gam(p^r)-\Gamma(p^{k+1})/\Gam(p^r)
&=A_k\cup\Big(\bigcup_{m=1}^{r-k}B_k^{(m)}\Big)\cup C_k,
\end{align*}
and 
\begin{align}
\#A_0^{(k,l)}&=\begin{cases}\disp\frac{1}{2}\varphi(l)p^{3r-k-2}(p^2-1)&(k<r),\\
\disp\frac{1}{2}\varphi(l)p^{2r-1}(p+1)&(k=r),\end{cases}&
\#A_k&=\frac{1}{2}p^{3r-3k-2}(p^2-1),\nt\\
\#B_k^{(m)}&=\begin{cases}p^{3r-3k-m-3}(p-1)^2(p+1)&(m<r-k),\\ 
p^{2r-2k-2}(p^2-1)&(m=r-k),\end{cases}\label{number}\\
\#B_k^{(m,+)}&=\#B_k^{(m,-)}=\frac{1}{2}\#B_k^{(m)},\nt\\
\#C_0^{(k,l)}&=\begin{cases}\disp\frac{1}{2}\varphi(l)p^{3r-k-2}(p-1)^2&(k<r),\\
\disp\frac{1}{2}\varphi(l)p^{2r-1}(p-1)&(k=r),\end{cases}&
\#C_k&=\frac{1}{2}p^{3r-3k-2}(p-1)^2,\nt
\end{align}
where $\varphi(l)$ is the Euler function given by 
$\varphi(l):=\#(\bZ/l\bZ)^*$. 
We also notice that the following relations hold among the sets 
$A_0^{(k,l)},A_k,B_k^{(m)},C_0^{(k,l)},C_k$ 
for $1\leq M\leq p$.
\begin{align}\label{str}
\{\gamma^M\divset\gamma\in A_0^{(k,l)}\}=&
\begin{cases}\Gamma(p^r)&(l|M,k=r),\\
A_k&(l|M,k\leq r-1),\\ 
A_0^{(k+1,l)}&(M=p,k\leq r-1),\\
A_0^{(k,l/\gcd{(M,l)})}&(\text{otherwise}),
\end{cases}\nt\\
\{\gamma^M\divset\gamma\in A_k\}=&
\begin{cases}\Gamma(p^r)&(M=p,k=r-1),\\
A_{k+1}&(M=p,k\leq r-2),\\ 
A_k&(\text{otherwise}),
\end{cases}\nt\\
\{\gamma^M\divset\gamma\in B_k^{(m,\pm)}\}=&
\begin{cases}\Gamma(p^r)&(M=p,k=r-1,m=1),\\
B_{k+1}^{(r-k-1)}&(M=p,k\leq r-2,m=r-k),\\ 
B_{k+1}^{(m,\pm)}&(M=p,k\leq r-2,m\leq r-k-1),\\ 
B_{k}^{(m,\pm)}&(\text{otherwise}),
\end{cases}
\end{align}
\begin{align}
\{\gamma^M\divset\gamma\in C_0^{(k,l)}\}=&
\begin{cases}\Gamma(p^r)&(l|M,k=r),\\
C_k&(l|M,m\leq r-1),\\ 
C_0^{(k+1,l)}&(M=p,k\leq r-1),\\
C_0^{(k,l/\gcd{(M,l)})}&(\text{otherwise}),
\end{cases}\nt\\
\{\gamma^M\divset\gamma\in C_k\}=&
\begin{cases}\Gamma(p^r)&(M=p,k=r-1),\\
C_{k+1}&(M=p,k\leq r-2),\\ 
C_k&(\text{otherwise}).
\end{cases}\nt
\end{align}

To calculate the type of the representatives of conjugacy classes of $\Xi$, 
the following lemma (when $r=1$, see also \cite{H}) is an inevitable step.
\begin{lem}\label{trg}
Let $\sigma(\gamma,\Gamma):=\tr\Big(\Ind_{\Gamma}^{\sz}1\Big)(\gamma)$. Then we have
\begin{align*}
\sigma(\gamma,\Gamma_0(p^r))=&\begin{cases}p^{r-1}(p+1)&(\gamma\in \Gamma(p^r)),\\
2p^k&(\gamma\in A_k),\\
2&(\gamma\in A_0^{(k,l)},\\
p^{[(r+k)/2]}&(\gamma\in B_k^{(r-k)}),\\
2p^{k+m/2}&(\gamma\in B_k^{(m,+)},\text{$m$ is even}),\\
0&(\text{otherwise}),
\end{cases} \\
\sigma(\gamma,\Gamma_1(p^r))=&\begin{cases}\disp\frac{1}{2}p^{2r-2}(p^2-1)&(\gamma\in\Gamma(p^r)),\\
\disp\frac{1}{2}p^{r+k-1}(p-1)&(\gamma\in B_k^{(r-k)}),\\
0&\text{(otherwise)}.\end{cases}
\end{align*}
\end{lem}

\begin{proof}
It is easy to see that the complete system of representatives of $\sz/\Gamma_0(p^r)$ 
can be chosen as
\begin{align}
\Bigg\{&\begin{pmatrix}1 &0\\ m& 1\end{pmatrix},\begin{pmatrix}lp &-1\\ 1& 0\end{pmatrix} 
\Biggdivset m\in\bZ/p^r\bZ, l\in\bZ/p^{r-1}\bZ\bigg\}.\label{csr}
\end{align}
Hence, we have
\begin{align*}
\sigma(\gamma,\Gamma_0(p^r))=&
\#\{m\in\bZ/p^r\bZ\divset \gamma_{12}m^2+(\gamma_{11}-\gamma_{22})m-\gamma_{21}\equiv0\bmod{p^r}\}\\
+&\#\{l\in\bZ/p^{r-1}\bZ\divset 
\gamma_{21}p^2l^2+p(\gamma_{11}-\gamma_{22})l-\gamma_{12}\equiv0\bmod{p^r}\}.
\end{align*}
Calculating the terms above for each element in the table of Lemma \ref{conjugacy}, 
we get the claims for $\Gamma=\Gamma_0(p^r)$.
Moreover, since $\Gamma_1(p^r)$ is a normal subgroup of $\Gamma_0(p^r)$, we have
\begin{align*}
\sigma(\gamma,\Gamma_1(p^r))&=\sum_{\begin{subarray}{c}g\in \CSR[\sz/\Gamma_0(p^r)]\\ 
g^{-1}\gamma g^{-1}\in \Gamma_0(p^r)\end{subarray}}
\tr\Big(\Ind_{\Gamma_1(p^r)}^{\Gamma_0(p^r)}1\Big)(g^{-1}\gamma g)\\
=&\sum_{\begin{subarray}{c}g\in \CSR[\sz/\Gamma_0(p^r)]\\ 
g^{-1}\gamma g^{-1}\in\Gamma_1(p^r)\end{subarray}}
[\Gamma_0(p^r):\Gamma_1(p^r)]\\
=&\frac{1}{2}p^{r-1}(p-1)\#\{g\in\CSR[\sz/\Gamma_0(p^r)]\bigdivset g^{-1}\gamma g\in\Gamma_1(p^r)\}.
\end{align*}
Hence the assertions for $\Gamma=\Gamma_1(p^r)$ follows from \eqref{csr}.
\end{proof}

By using \eqref{str} and the lemma above, we determine the type of each element 
of $\Xi$ as follows.
\begin{lem}\label{type}
We have
\begin{align*}
\lambda_0^{p^r}(1)=&\text{type of $\gamma\in\Gamma(p^r)$ in $\Gamma_0(p^r)$},\\
\lambda_0^{p^r}(lp^{k})=&
\begin{cases}\text{type of $\gamma\in A_0^{(r-k,l)}$ in $\Gamma_0(p^r)$} &(l|(p-1)/2,l>1),\\
\text{type of $\gamma\in C_0^{(r-k,l)}$ in $\Gamma_0(p^r)$} &(l|(p+1)/2,l>1),\end{cases}\\
\lambda_0^{p^r}(p^{k},A)=&\text{type of $\gamma\in A_{r-k}$ in $\Gamma_0(p^r)$},\\
\lambda_0^{p^r}(p^{k},B^{(m,\pm)})
=&\text{type of $\gamma\in B_{r-k}^{(m,\pm)}$ in $\Gamma_0(p^r)$},\\
\lambda_0^{p^r}(p^{k},C)=&\text{type of $\gamma\in C_{r-k}$ in $\Gamma_0(p^r)$},
\end{align*}
and 
\begin{align*}
\lambda_1^{p^r}(1)=&\text{type of $\gamma\in\Gamma(p^r)$ in $\Gamma_1(p^r)$},\\
\lambda_1^{p^r}(lp^{k})=&
\begin{cases}\text{type of $\gamma\in A_0^{(r-k,l)}$ in $\Gamma_1(p^r)$} &(l|(p-1)/2,l>1),\\
\text{type of $\gamma\in C_0^{(r-k,l)}$ in $\Gamma_1(p^r)$} &(l|(p+1)/2,l>1),\end{cases}\\
\lambda_1^{p^r}(p^{k})=&\text{type of $\gamma\in A_{r-k}\cup  C_{r-k}$ in $\Gamma_1(p^r)$},\\
\lambda_1^{p^r}(p^{k},B^{(m)})=&\text{type of $\gamma\in B_{r-k}^{(m)}$ in $\Gamma_1(p^r)$}.
\end{align*}
\end{lem}
\begin{proof}
By Lemma \ref{prop}, if $\gamma$ is $(1^{l_1}2^{l_2}\cdots n^{l_n})$-type in $\tilde{\Gamma}$, 
then we have $\tr{\sigma(\gamma^m)}=\sum_{j|m}jl_j$. 
Hence, $l_j$'s are calculated recursively by 
\begin{align}\label{ind2}
ml_m=\sum_{j|m}\mu(m/j)\tr\sigma(\gamma^{m/j}), 
\end{align}
where $mu$ is the M\"{o}bius function.
Hence, using \eqref{str}, \eqref{ind2} and Lemma \ref{trg}, 
we obtain the desired results recursively. 
This proves the lemma.
\end{proof}

Applying Lemma \ref{type} and \eqref{number} to Corollary \ref{cor}, we can calculate the densities 
for $N=p^r$ when $p$ is an odd prime.
The densities for $N=2^r$ can be calculated similarly.
For a general integer $N>1$, 
by the help of the following proposition, 
we determine $\pi_{\tilde{\Gamma}\upar \sz}^{\lambda}(x)$ recursively.

\begin{prop}\label{last}
Let $\Gamma_1$ and $\Gamma_2$ be subgroups of $\Gamma$ of finite index, and
$\Gamma'_1$ and $\Gamma'_2$ the maximal normal subgroups of $\Gamma$ contained 
in $\Gamma_1$ and $\Gamma_2$ respectively. 
Put $\Xi_1:=\Gamma/\Gamma'_1$, $\Xi_2:=\Gamma/\Gamma'_2$, 
$n_1:=[\Gamma:\Gamma_1]$ and $n_2:=[\Gamma:\Gamma_2]$.
Assume that $\Gamma'_1$ and $\Gamma'_2$ are relatively prime in $\Gamma$, i.e.
the index of $\Gamma'_1\cap\Gamma'_2$ in $\Gamma$ is finite 
and $\Gamma'_1\Gamma'_2=\Gamma$ 
(note that $\Gamma_1$ and $\Gamma_2$ are also relatively prime in $\Gamma$). 
Then we have
\begin{align*}
\lim_{x\tinf}\mu_{\Gamma_1\cap\Gamma_2\upar \Gamma}^{\lambda}(x)
=\sum_{\begin{subarray}{c}\lambda_1\vdash n_1,\lambda_2\vdash n_2
\\ \lambda=\lambda_1\otimes\lambda_2\end{subarray}}
\lim_{x\tinf}\mu_{\Gamma_1\upar \Gamma}^{\lambda_1}(x)
\lim_{x\tinf}\mu_{\Gamma_2\upar \Gamma}^{\lambda_2}(x),
\end{align*}
where $\lambda_1\otimes\lambda_2:=
(m_1l_1,\cdots,m_1l_{k_2},m_2l_1,\cdots,m_{k_1}l_{k_2})\vdash n_1n_2$ 
for $\lambda_1=(m_1,\cdots,m_{k_1})\vdash n_1$ 
and $\lambda_2=(l_1,\cdots,l_{k_2})\vdash n_2$. 
For a partition $\lambda$ which can not be expressed as $\lambda_1\otimes\lambda_2$, 
we have $\mu_{\Gamma_1\cap\Gamma_2\upar \Gamma}^{\lambda}(x)=0$.
\end{prop}

\begin{rem}
\rm{It is clear that, if $\Gamma'_1$ and $\Gamma'_2$ are relatively prime in $\Gamma$, 
then $\Gamma_1$ and $\Gamma_2$ are also relatively prime in $\Gamma$. 
The converse is true in the cases where $\Gamma=\sz$ and $\Gamma_1,\Gamma_2$ 
are congruence subgroups of $\Gamma$ with relatively prime levels. 
In general, however, this is not true; 
in fact, if we take $\Gamma=\Gamma_0(p)$, $\Gamma_1=\Gamma_1(p)$ and}
\begin{align*}
\Gamma_2=\bigg\{\gamma\in \sz\Bigdivset
\gamma\equiv\begin{pmatrix}\delta  &0\\ 0& \delta^{-1}\end{pmatrix}\bmod{p},
\delta\in(\bZ/p\bZ)^{*}\bigg\},
\end{align*}
\rm{we have $\Gamma_1\Gamma_2=\Gamma$ but $\Gamma'_1\Gamma'_2\not=\Gamma$ 
because $\Gamma'_1=\Gamma_1(p)$ and $\Gamma'_2=\Gamma(p)$.}
\end{rem}

To show Proposition \ref{last}, we prepare the following lemmas.
\begin{lem}\label{lem11}
Let $\Gamma_1$ and $\Gamma_2$ be relatively prime subgroups of $\Gamma$. 
Then $[\Gamma:\Gamma_1\cap\Gamma_2]=[\Gamma:\Gamma_1][\Gamma:\Gamma_2]$ holds, 
and $\CSR[\Gamma_2/(\Gamma_1\cap\Gamma_2)]$ gives 
a complete system of representatives of $\Gamma/\Gamma_1$.
\end{lem}
\begin{proof}
Let $\CSR[\Gamma_2/(\Gamma_1\cap\Gamma_2)]=\{a_1,\cdots,a_k\}$.
It is easy to see that 
\begin{align*}
\Gamma_2=&\bigcup_{i=1}^{k}a_i(\Gamma_1\cap\Gamma_2)=\bigcup_{i=1}^{k}(a_i\Gamma_1\cap\Gamma_2)
=\big(\bigcup_{i=1}^{k}a_i\Gamma_1\big)\cap\Gamma_2.
\end{align*}
Hence, we have $\bigcup_{i=1}^{k}a_i\Gamma_1\supset\Gamma_2$.
Since $\Gamma_1$ and $\Gamma_2$ are relatively prime, 
and $\cup_{i=1}^{k}a_i\Gamma_1$ contains both $\Gamma_1$ and $\Gamma_2$, 
we have $\bigcup_{i=1}^{k}a_i\Gamma_1=\Gamma$. 
Hence $\CSR[\Gamma/\Gamma_1]$ can be chosen as a subset of $\CSR[\Gamma_2/(\Gamma_1\cap\Gamma_2)]$.
Now, we choose $\CSR[\Gamma/\Gamma_1]:=\{b_1,\cdots,b_l\}
\subset\CSR[\Gamma_2/(\Gamma_1\cap\Gamma_2)]\subset\Gamma_2$. 
Then, it is easy to see that 
\begin{align*}
\Gamma_2=\Gamma\cap\Gamma_2=&\Big(\bigcup_{j=1}^{l}b_j\Gamma_1\Big)\cap\Gamma_2
=\bigcup_{j=1}^{l}(b_j\Gamma_1\cap\Gamma_2)=\bigcup_{j=1}^{l}b_j(\Gamma_1\cap\Gamma_2).
\end{align*} 
Therefore, we conclude that $\CSR[\Gamma/\Gamma_1]=\CSR[\Gamma_2/(\Gamma_1\cap\Gamma_2)]$. 
\end{proof}

\begin{lem}\label{lem12}
Let $\Gamma_1$ and $\Gamma_2$ be relatively prime subgroups of $\Gamma$. 
If $\gamma\in\Prim(\Gamma)$ is $\lambda_1$-type in $\Gamma_1$ 
and is also $\lambda_2$-type in $\Gamma_2$, 
then $\gamma$ is $\lambda_1\otimes\lambda_2$-type in $\Gamma_1\cap\Gamma_2$. 
Furthermore, if $\gamma$ is $\lambda$-type in $\Gamma_1\cap\Gamma_2$, 
then there exist $\lambda_1\vdash n_1$, $\lambda_2\vdash n_2$ 
such that $\lambda=\lambda_1\otimes\lambda_2$, 
and $\gamma$ is simultaneously $\lambda_1$-type in $\Gamma_1$ and $\lambda_2$-type in $\Gamma_2$.
\end{lem}
\begin{proof} 
By Lemma \ref{lem11}, it is easy to see that  
\begin{align*}
\Ind_{\Gamma_1\cap\Gamma_2}^{\Gamma}1
=\Ind_{\Gamma_2}^{\Gamma}\big(\Ind_{\Gamma_1\cap\Gamma_2}^{\Gamma_2}1\big)
=\Ind_{\Gamma_2}^{\Gamma}\big(\Ind_{\Gamma_1}^{\Gamma}1\big|_{\Gamma_2}\big).
\end{align*}
Hence, if we put $\CSR[\Gamma/\Gamma_2]=\{a_1,\cdots,a_{n_2}\}$, 
we have 
\begin{align*}
\tr\big(\Ind_{\Gamma_1\cap\Gamma_2}^{\Gamma}1\big)(\gamma)
&=\tr\Big(\Ind_{\Gamma_2}^{\Gamma}\big(\Ind_{\Gamma_1}^{\Gamma}1\big|_{\Gamma_2}\big)\Big)(\gamma)\\
&=\sum_{\begin{subarray}{c}1\leq i\leq n_2\\ a_i^{-1}\gamma a_i\in\Gamma_2\end{subarray}}
\tr\big(\Ind_{\Gamma_1}^{\Gamma}1\big|_{\Gamma_2}\big)(a_i^{-1}\gamma a_i)\\
&=\sum_{\begin{subarray}{c}1\leq i\leq n_2\\ a_i^{-1}\gamma a_i\in\Gamma_2\end{subarray}}
\tr\big(\Ind_{\Gamma_1}^{\Gamma}1\big)(\gamma)
=\tr\big(\Ind_{\Gamma_2}^{\Gamma}1\big)(\gamma)\times
\tr\big(\Ind_{\Gamma_1}^{\Gamma}1\big)(\gamma).
\end{align*}
Similarly, for $k\geq1$, we also have 
\begin{align*}
\tr\big(\Ind_{\Gamma_1\cap\Gam_2}^{\Gamma}1\big)(\gamma^k)
=\tr\big(\Ind_{\Gamma_2}^{\Gamma}1\big)(\gamma^k)\times
\tr\big(\Ind_{\Gamma_1}^{\Gamma}1\big)(\gamma^k).
\end{align*}
Hence, according to \eqref{ind2}, we see that the type of 
$\big(\Ind_{\Gamma_1\cap\Gamma_2}^{\Gamma}1\big)(\gamma)$ coincides that of 
$\big(\Ind_{\Gamma_1}^{\Gamma}1\big)(\gamma)\otimes\big(\Ind_{\Gamma_2}^{\Gamma}1\big)(\gamma)$. 
This completes the proof of the lemma.
\end{proof}

\noindent{\bf Proof of Proposition \ref{last}.}
Since $\Gamma'_1$ and $\Gamma'_2$ are relatively prime, 
according to Lemma \ref{lem11}, 
we can choose $\CSR[\Gamma/\Gamma'_1]$ and $\CSR[\Gamma/\Gamma'_2]$  
as subsets of $\Gamma'_2$ and $\Gamma'_1$ respectively. 
Now, we put $\CSR[\Gamma/\Gamma'_1]=\{a_1,\cdots,a_{|\Xi_1|}\}\subset\Gamma'_2$ 
and $\CSR[\Gamma/\Gamma'_2]:=\{b_1,\cdots,b_{|\Xi_2|}\}\subset\Gamma'_1$. 
Then $\CSR[\Gamma/(\Gamma'_1\cap\Gamma'_2)]$ can be chosen as 
$\{a_i b_j\}_{1\leq i\leq |\Xi_1|,1\leq j\leq |\Xi_2|}$. 
For $\lambda_1\vdash n_1$ and $\lambda_2\vdash n_2$, 
we denote by $a'_1,\cdots,a'_{k_1}$ (resp. $b'_1,\cdots,b'_{k_2}$) 
the elements of $\CSR[\Gamma/\Gamma'_1]$ (resp. $\CSR[\Gamma/\Gamma'_2]$
which are $\lambda_1$-type in $\Gamma_1$ (resp. $\lambda_2$-type in $\Gamma_2$). 
It is easy to see that $a'_i b'_j$ ($1\leq i\leq k_1$, $1\leq j\leq k_2$) 
is $\lambda_1$-type in $\Gamma_1$ 
and is $\lambda_2$-type in $\Gamma_2$. 
Hence, it follows from Theorem \ref{thm1} that
\begin{align*}
&\#\{\gamma\in\CSR[\Gamma/(\Gamma'_1\cap\Gamma'_2)]\divset
\text{$\gamma$ is $\lambda_1$-type in $\Gamma_1$ 
and is $\lambda_2$-type in $\Gamma_2$}\}\\
=&\#\{\gamma\in\CSR[\Gamma/\Gamma'_1]\divset\text{$\gamma$ is $\lambda_1$-type in $\Gamma_1$}\}
\#\{\gamma\in\CSR[\Gamma/\Gamma'_2]\divset\text{$\gamma$ is $\lambda_2$-type in $\Gamma_2$}\}\\
=&|\Xi_1|\lim_{x\tinf}\mu_{\Gamma_1}^{\lambda_1}(x)
\times|\Xi_2|\lim_{x\tinf}\mu_{\Gamma_2}^{\lambda_2}(x).
\end{align*}
Therefore, by Lemma \ref{lem12}, we have the proposition. 
\qed

\section{Examples for congruence subgroups}

\noi {\bf The case of $\tilde{\Gamma}=\Gamma_0(3)$.}\\
In this case, $\Gamma'=\Gamma(3)$, $\Xi=\mathrm{SL}_2(\bZ/3\bZ)/\{\pm\Id\}$, $|\Xi|=12$ and $n=1$.
We have 
\begin{align*}
\lim_{x\tinf}\mu_{\Gamma_0(3)\upar \sz}^{\lambda^3_0(1)}(x)=&\frac{1}{12},\quad
\lim_{x\tinf}\mu_{\Gamma_0(3)\upar \sz}^{\lambda^3_0(2)}(x)=\frac{1}{4},\\
\lim_{x\tinf}\mu_{\Gamma_0(3)\upar \sz}^{\lambda^3_0(3)}(x)=&\frac{2}{3},
\end{align*}
where 
\begin{align*}
\lambda^3_0(1)=&(1^4)\quad\big(=\miniyng(1,1,1,1)\quad\big),\quad
\lambda^3_0(2)=(2^2)\quad\big(=\miniyng(2,2)\quad\big),\\
\lambda^3_0(3)=&(3,1)\quad\big(=\miniyng(3,1)\quad\big).
\end{align*}

\noi {\bf The case of $\tilde{\Gamma}=\Gamma_0(5)$.}\\
In this case, $\Gamma'=\Gamma(5)$, $\Xi=\mathrm{SL}_2(\bZ/5\bZ)/\{\pm\Id\}$, $|\Xi|=60$ and $n=6$.
We have 
\begin{align*}
&\lim_{x\tinf}\mu_{\Gamma_0(5)\upar \sz}^{\lambda^5_0(1)}(x)=\frac{1}{60},&&
\lim_{x\tinf}\mu_{\Gamma_0(5)\upar \sz}^{\lambda^5_0(2)}(x)=\frac{1}{4},\\
&\lim_{x\tinf}\mu_{\Gamma_0(5)\upar \sz}^{\lambda^5_0(3)}(x)=\frac{1}{3},&&
\lim_{x\tinf}\mu_{\Gamma_0(5)\upar \sz}^{\lambda^5_0(5)}(x)=\frac{2}{5},
\end{align*}
where 
\begin{align*}
&\lambda^5_0(1)=(1^6)\quad\big(=\miniyng(1,1,1,1,1,1)\quad\big),&&
\lambda^5_0(2)=(2^2,1^2)\quad\big(=\miniyng(2,2,1,1)\quad\big),\\
&\lambda^5_0(3)=(3^2)\quad\big(=\miniyng(3,3)\quad\big),&&
\lambda^5_0(5)=(5,1)\quad\big(=\miniyng(5,1)\quad\big).
\end{align*}

\noi {\bf The case of $\tilde{\Gamma}=\Gamma_0(5^2)$.}\\
In this case, $\Gamma'=\Gamma(5^2)$, $\Xi=\mathrm{SL}_2(\bZ/5^2\bZ)/\{\pm\Id\}$, 
$|\Xi|=7500$ and $n=30$.
We have 
\begin{align*}
\lim_{x\tinf}\mu_{\Gamma_0(25)\upar \sz}^{\lambda^{25}_0(1)}(x)=&\frac{1}{7500},&
\lim_{x\tinf}\mu_{\Gamma_0(25)\upar \sz}^{\lambda^{25}_0(2)}(x)=&\frac{1}{20},\\
\lim_{x\tinf}\mu_{\Gamma_0(25)\upar \sz}^{\lambda^{25}_0(3)}(x)=&\frac{1}{15},&
\lim_{x\tinf}\mu_{\Gamma_0(25)\upar \sz}^{\lambda^{25}_0(5,A)}(x)=&\frac{1}{125},\\
\lim_{x\tinf}\mu_{\Gamma_0(25)\upar \sz}^{\lambda^{25}_0(5,B^{(1)})}(x)=&\frac{2}{625},&
\lim_{x\tinf}\mu_{\Gamma_0(25)\upar \sz}^{\lambda^{25}_0(5,C)}(x)=&\frac{2}{375},\\
\lim_{x\tinf}\mu_{\Gamma_0(25)\upar \sz}^{\lambda^{25}_0(10)}(x)=&\frac{1}{5},&
\lim_{x\tinf}\mu_{\Gamma_0(25)\upar \sz}^{\lambda^{25}_0(15)}(x)=&\frac{4}{15},\\
\lim_{x\tinf}\mu_{\Gamma_0(25)\upar \sz}^{\lambda^{25}_0(25,B^{(1)})}(x)=&\frac{8}{25},&
\lim_{x\tinf}\mu_{\Gamma_0(25)\upar \sz}^{\lambda^{25}_0(25,B^{(2)})}(x)=&\frac{2}{25},
\end{align*}
where 
\begin{align*}
\lambda^{25}_0(1)=&(1^{30}),&
\lambda^{25}_0(2)=&(2^{14},1^2),\\
\lambda^{25}_0(3)=&(3^{10}),&
\lambda^{25}_0(5,A)=&(5^4,1^{10}),\\
\lambda^{25}_0(5,B^{(1)})=&(5^5,1^5),&
\lambda^{25}_0(5,C)=&(5^{6}),\\
\lambda^{25}_0(10)=&(10^2,2^4,1^2),&
\lambda^{25}_0(15)=&(15^2),\\
\lambda^{25}_0(25,B^{(1)})=&(25,1^5),&
\lambda^{25}_0(25,B^{(2)})=&(25,5).
\end{align*}

\noi {\bf The case of $\tilde{\Gamma}=\Gamma_0(3\times5^2)$.}\\
Since $\Gamma_0(3\times5^2)=\Gamma_0(3)\cap\Gamma_0(5^2)$ and 
$\Gamma'=\Gamma(3)\cap\Gamma(5^2)$, 
by employing Proposition \ref{last}, we have the following results ($|\Xi|=90000,n=120$). 
\begin{align*}
\lim_{x\tinf}\mu_{\Gamma_0(75)\upar \sz}^{(1^{120})}(x)=&\frac{1}{90000},&
\lim_{x\tinf}\mu_{\Gamma_0(75)\upar \sz}^{(2^{56},1^8)}(x)=&\frac{1}{240},\\
\lim_{x\tinf}\mu_{\Gamma_0(75)\upar \sz}^{(10^8,2^{16},1^8)}(x)=&\frac{1}{60},&
\lim_{x\tinf}\mu_{\Gamma_0(75)\upar \sz}^{(3^{40})}(x)=&\frac{1}{180},\\
\lim_{x\tinf}\mu_{\Gamma_0(75)\upar \sz}^{(15^{8})}(x)=&\frac{1}{45},&
\lim_{x\tinf}\mu_{\Gamma_0(75)\upar \sz}^{(5^{16},1^{40})}(x)=&\frac{1}{1500},\\
\lim_{x\tinf}\mu_{\Gamma_0(75)\upar \sz}^{(5^{20},1^{20})}(x)=&\frac{1}{3750},&
\lim_{x\tinf}\mu_{\Gamma_0(75)\upar \sz}^{(25^4,1^{20})}(x)=&\frac{2}{75},\\
\lim_{x\tinf}\mu_{\Gamma_0(75)\upar \sz}^{(25^4,5^{2})}(x)=&\frac{1}{150},&
\lim_{x\tinf}\mu_{\Gamma_0(75)\upar \sz}^{(5^{24})}(x)=&\frac{1}{2250},\\
\lim_{x\tinf}\mu_{\Gamma_0(75)\upar \sz}^{(2^{60})}(x)=&\frac{1}{30000},&
\lim_{x\tinf}\mu_{\Gamma_0(75)\upar \sz}^{(4^{28},2^4)}(x)=&\frac{1}{80},\\
\lim_{x\tinf}\mu_{\Gamma_0(75)\upar \sz}^{(20^4,4^8,2^4)}(x)=&\frac{1}{20},&
\lim_{x\tinf}\mu_{\Gamma_0(75)\upar \sz}^{(6^{20})}(x)=&\frac{1}{60},\\
\lim_{x\tinf}\mu_{\Gamma_0(75)\upar \sz}^{(30^4)}(x)=&\frac{1}{15},&
\lim_{x\tinf}\mu_{\Gamma_0(75)\upar \sz}^{(10^8,2^{10})}(x)=&\frac{1}{500},\\
\lim_{x\tinf}\mu_{\Gamma_0(75)\upar \sz}^{(10^{10},2^{10})}(x)=&\frac{1}{1250},&
\lim_{x\tinf}\mu_{\Gamma_0(75)\upar \sz}^{(50^2,2^{10})}(x)=&\frac{2}{25},\\
\lim_{x\tinf}\mu_{\Gamma_0(75)\upar \sz}^{(50^{2},10^2)}(x)=&\frac{1}{50},&
\lim_{x\tinf}\mu_{\Gamma_0(75)\upar \sz}^{(10^{12})}(x)=&\frac{1}{750},\\
\lim_{x\tinf}\mu_{\Gamma_0(75)\upar \sz}^{(3^{30},1^{30})}(x)=&\frac{1}{11250},&
\lim_{x\tinf}\mu_{\Gamma_0(75)\upar \sz}^{(6^{14},3^2,2^{14},1^2)}(x)=&\frac{1}{30},\\
\lim_{x\tinf}\mu_{\Gamma_0(75)\upar \sz}^{(30^2,10^2,6^4,3^2,2^4,1^2)}(x)=&\frac{2}{15},&
\lim_{x\tinf}\mu_{\Gamma_0(75)\upar \sz}^{(9^{10},3^{10})}(x)=&\frac{2}{45},\\
\lim_{x\tinf}\mu_{\Gamma_0(75)\upar \sz}^{(45^2,15^2)}(x)=&\frac{8}{45},&
\lim_{x\tinf}\mu_{\Gamma_0(75)\upar \sz}^{(15^4,5^5,3^{10},1^{10})}(x)=&\frac{2}{375},\\
\lim_{x\tinf}\mu_{\Gamma_0(75)\upar \sz}^{(15^5,5^5,3^5,1^5)}(x)=&\frac{4}{1875},&
\lim_{x\tinf}\mu_{\Gamma_0(75)\upar \sz}^{(75,25,3^5,1^5)}(x)=&\frac{16}{75},\\
\lim_{x\tinf}\mu_{\Gamma_0(75)\upar \sz}^{(75,25,15,5)}(x)=&\frac{4}{75},&
\lim_{x\tinf}\mu_{\Gamma_0(75)\upar \sz}^{(15^6,5^6)}(x)=&\frac{4}{1125}.
\end{align*}

\begin{rem}
\rm{In the case $\tilde{\Gamma}=\Gamma_0(3\times5^2)$, 
$\lambda^3_0(m_1)\otimes\lambda^{5^2}_0(l_1)\not=\lambda^3_0(m_2)\otimes\lambda_0^{5^2}(l_2)$ 
holds if $(m_1,l_1)\not=(m_2,l_2)$. 
For a general pair of relatively prime $\Gamma_1$ and $\Gamma_2$, 
we do not know whether there are partition $\lambda_1,\lambda'_1\vdash n_1$ 
and $\lambda_2,\lambda'_2\vdash n_2$ 
($\lambda_1\not=\lambda'_1$, $\lambda_2\not=\lambda'_2$) 
such that $\lambda_1\otimes\lambda_2=\lambda'_1\otimes\lambda'_2$ provided 
$\mu_{\Gamma_1\upar \Gamma}^{\lambda_1}(x)$, $\mu_{\Gamma_1\upar \Gamma}^{\lambda'_1}(x)$, 
$\mu_{\Gamma_2\upar \Gamma}^{\lambda_2}(x)$ and $\mu_{\Gamma_2\upar \Gamma}^{\lambda'_2}(x)$ 
having non-zero densities.}
\end{rem}

\begin{rem}
\rm{All elements $\gamma$ of $\Gamma$ are $(\cdots,M(\gamma)^l)$-type in $\tilde{\Gamma}$ ($l>0$) 
at the examples discussed in this section 
(see Theorem \ref{thm1} for the definition of $M(\gamma)$).
In general, we may expect that there is some $\gamma$ whose type is of the form $(\cdots,M)$ 
for $M<M(\gamma)$, but we have never found such examples unfortunately. 
It is interesting to study the density  
$\hat{\mu}_{\tilde{\Gamma}\upar\Gamma}(x):=\hat{\pi}_{\tilde{\Gamma}\upar\Gamma}(x)/\pi_{\Gamma}(x)$ 
as $x\tinf$, 
where $\hat{\pi}_{\tilde{\Gamma}\upar\Gamma}(x):=\{\gamma\in\Prim(\Gamma)\divset N(\gamma)<x,
\text{$\gamma$ is $(\cdots,M)$-type in $\tilde{\Gamma}$ for $\exists M<M(\gamma)$}\}$.}
\end{rem}

\section{A Remark on Selberg's zeta functions}
Let $\zeta_{\Gamma}(s)$ be a Selberg zeta function of $\Gamma$ defined by
\begin{align*}
\zeta_{\Gamma}(s):=\prod_{\gamma\in\Prim(\Gamma)}(1-N(\gamma)^{-s})^{-1}\quad \Re{s}>1.
\end{align*} 
For $\tilde{\Gamma}\subset \Gamma$ and $\lambda\vdash n$, we define 
a Selberg type zeta function attached to this data by
\begin{align*}
\zeta_{\tilde{\Gamma}\upar\Gamma}^{\lambda}(s):=\prod_{\begin{subarray}{c}\gamma\in\Prim(\Gamma)\\
\text{$\gamma$ is $\lambda$-type in $\tilde{\Gamma}$}\end{subarray}}(1-N(\gamma)^{-s})^{-1}
\quad \Re{s}>1.
\end{align*}
By using Venkov-Zograf's formula \cite{VZ}, we have 
\begin{align}
\zeta_{\tilde{\Gamma}}(s)=&\zeta_{\Gamma}(s,\sigma)\nt\\
=&\prod_{\gamma\in\Prim(\Gamma)}\det{(\Id-\sigma(\gamma)N(\gamma)^{-s})}^{-1}\nt\\
=&\prod_{\lambda\vdash n}\prod_{\begin{subarray}{c}\gamma\in\Prim(\Gamma)\\
\text{$\gamma$ is $\lambda$-type in $\tilde{\Gamma}$}\end{subarray}}
\det{(\Id-\sigma(\gamma)N(\gamma)^{-s})}^{-1}\nt\\
=&\prod_{\lambda=(m_1,m_2,\cdots,m_k)\vdash n}\zeta_{\tilde{\Gamma}\upar\Gamma}^{\lambda}(m_1s)
\cdots\zeta_{\tilde{\Gamma}\upar\Gamma}^{\lambda}(m_ks),\label{venkov}
\end{align}
where $\sigma:=\Ind_{\tilde{\Gam}}^{\Gam}1$.
Although detailed studies of $\zeta_{\tilde{\Gamma}\upar\Gamma}^{\lambda}(s)$ remain in the future, 
in this section, we study $\zeta_{\tilde{\Gamma}\upar\Gamma}^{\lambda}(s)$ for  
the particular cases where $\Gamma=\sz$, $\tilde{\Gamma}=\Gamma_1(p)$ and $\Gamma(p)$ as follows.

In these cases, from the discussions in Section 4, it is easy to see that
\begin{align*}
\text{$\gamma$ is $\lambda_1^p(m)$-type in $\Gamma_1(p)$}\Leftrightarrow
\text{$\gamma$ is $\lambda^p(m)$-type in $\Gamma(p)$}\Leftrightarrow M(\gamma)=m.
\end{align*}
Then we have 
\begin{align*}
\zeta_{\Gamma_1(p)\upar \sz}^{\lambda_1^p(m)}(s)=\zeta_{\Gamma(p)\upar \sz}^{\lambda^p(m)}(s)
=\prod_{\begin{subarray}{c}\gamma\in\Prim(\Gamma)\\
M(\gamma)=m\end{subarray}}(1-N(\gamma)^{-s})^{-1}.
\end{align*}
For simplicity, we denote this function by $\zeta_{\sz}^{(p,m)}(s)$. 
The functions $\zeta_{\sz}^{(p,m)}(s)$ have the following properties.
\begin{prop}\label{selberg}
Let $p$ be an odd prime. Then we have
\begin{align}
\bigg\{\frac{\big(\zeta_{\sz}^{(p,p)}(s)\big)^p}{\zeta_{\sz}^{(p,p)}(ps)}\bigg\}^{\frac{p-1}{2}}
=\frac{\big(\zeta_{\Gamma_1(p)}(s)\big)^p}{\zeta_{\Gamma(p)}(s)}.\label{ratio}
\end{align}
Furthermore, $\big(\zeta_{\sz}^{(p,p)}(s)\big)^{p^{r}(p-1)/2}$ can be analytically continued 
to $\Re{s}>1/p^{r}$ as a meromorphic function and 
$\zeta_{\sz}^{(p,p)}(s)$ has infinitely many singular points near $s=0$.
\end{prop}

\begin{proof}
By using \eqref{venkov} and the results in Section 4, we have
\begin{align*}
\zeta_{\Gamma(p)}(s)=&\Big(\zeta_{\sz}^{(p,1)}(s)\Big)^{-\frac{1}{2}p(p^2-1)}
\Big(\zeta_{\sz}^{(p,p)}(ps)\Big)^{-\frac{1}{2}(p^2-1)}
\prod_{m|\frac{p\pm1}{2},m>1}\Big(\zeta_{\sz}^{(p,m)}(ms)\Big)^{-\frac{p(p^2-1)}{2m}},\\
\zeta_{\Gamma_1(p)}(s)=&\Big(\zeta_{\sz}^{(p,1)}(s)\Big)^{-\frac{1}{2}(p^2-1)}
\Big(\zeta_{\sz}^{(p,p)}(s)\zeta_{\sz}^{(p,p)}(ps)\Big)^{-\frac{1}{2}(p-1)}
\prod_{m|\frac{p\pm1}{2},m>1}\Big(\zeta_{\sz}^{(p,m)}(ms)\Big)^{-\frac{(p^2-1)}{2m}}.
\end{align*}
Hence the formula \eqref{ratio} follows immediately. 

Since $\zeta_{\Gamma_1(p)}(s)$ and $\zeta_{\Gamma(p)}(s)$ are meromorphic in the whole $\bC$  
and $\zeta_{\sz}^{(p,p)}(ps)$ is non-zero and holomorphic in $\Re{s}>1/p$, 
we see that $\big(\zeta_{\sz}^{(p,p)}(s)\big)^{p(p-1)/2}$ is analytically continued to $\Re{s}>1/p$ 
as a meromorphic function. 
Now, we take $p^{r-1}$-powers of the both hand sides of \eqref{ratio}. Then we have 
\begin{align}
\big(\zeta_{\sz}^{(p,p)}(s)\big)^{p^{r}(p-1)/2}
=&\big(\zeta_{\sz}^{(p,p)}(ps)\big)^{p^{r-1}(p-1)/2}\bigg\{\frac{\big(\zeta_{\Gamma_1(p)}(s)\big)^p}
{\zeta_{\Gamma(p)}(s)}\bigg\}^{p^{r-1}}\nt\\
=&\big(\zeta_{\sz}^{(p,p)}(p^2s)\big)^{p^{r-2}(p-1)/2}
\bigg\{\frac{\big(\zeta_{\Gamma_1(p)}(ps)\big)^p}{\zeta_{\Gamma(p)}(ps)}\bigg\}^{p^{r-2}}
\bigg\{\frac{\big(\zeta_{\Gamma_1(p)}(s)\big)^p}{\zeta_{\Gamma(p)}(s)}\bigg\}^{p^{r-1}}\nt\\
=&\cdots\nt\\
=&\big(\zeta_{\sz}^{(p,p)}(p^{r}s)\big)^{(p-1)/2}
\prod_{k=1}^{r}\bigg\{\frac{\big(\zeta_{\Gamma_1(p)}(p^{k}s)\big)^p}{\zeta_{\Gamma(p)}(p^{k}s)}
\bigg\}^{p^{r-k}}.
\label{ratio1}
\end{align}
Since $\zeta_{\sz}^{(p,p)}(p^{r}s)$ is non-zero holomorphic in $\Re{s}>1/p^r$, 
we can obtain the meromorphic continuation of $\big(\zeta_{\sz}^{(p,p)}(s)\big)^{p^{r}(p-1)/2}$ 
to the half plane $\Re{s}>1/p^{r}$. 

Both the functions $\zeta_{\Gamma_1(p)}(s)$ and $\zeta_{\Gamma(p)}(s)$ have simple poles 
at $s=1$ (see \cite{He}).
Hence, $\big(\zeta_{\sz}^{(p,p)}(s)\big)^{p}$ has a double pole at $s=1$.
Thus $\zeta_{\sz}^{(p,p)}(ps)$ has a branch point at $s=1/p$. 
Since $\zeta_{\Gamma_1(p)}(s)$ and $\zeta_{\Gamma(p)}(s)$ are meromorphic at $s=1/p$, 
by \eqref{ratio}, $\zeta_{\sz}^{(p,p)}(s)$ should have a branch point at $s=1/p$.
Then $\zeta_{\sz}^{(p,p)}(ps)$ has a branch point at $s=1/p^2$. 
Successively, we see that $\zeta_{\sz}^{(p,p)}(s)$ has branch points 
at $s=1,1/p,1/p^2,1/p^3,\cdots$. 
This shows that the series of the branch points $\{1,1/p,1/p^2,1/p^3,\cdots\}$ 
has an accumulation point at $0$.
Similarly, it is easy to see that $\zeta_{\sz}^{(p,p)}(s)$ has branch points 
at $s=1/2+ir_j,(1/2+ir_j)/p,(1/2+ir_j)/p^2,\cdots$,  
where $1/4+r_j^2$ is the $j$-th non-trivial eigenvalue of the Laplacian on $X_{\Gamma(p)}$ 
but not the spectrum of $X_{\Gamma_1(p)}$.
Hence $\zeta_{\sz}^{(p,p)}(s)$ has infinitely many branch points near $s=0$.
\end{proof}

\begin{rem}
\rm{In Proposition \ref{selberg}, we obtain the analytic continuation 
of $\zeta_{\sz}^{(p,p)}(s)$ to the half plane $\Re{s}>0$.
We do not know, however, whether $\zeta_{\sz}^{(p,p)}(s)$ can be analytically continued 
to a region contained in $\Re{s}\leq0$ 
or has a natural boundary $\Re{s}=0$.}
\end{rem}

\noindent{\bf Acknowledgement.} 
The authors would like to thank the referee for his/her useful comments.

\begin{flushleft}
\textsc{Graduate School of Mathematics\\ Kyushu University\\  
6-10-1, Hakozaki, Fukuoka, 812-8581\\ JAPAN}\\ 
\textit{E-mail address}: hasimoto@math.kyushu-u.ac.jp\\ 
\end{flushleft}
\begin{flushleft}
\textsc{Faculty of Mathematics\\ Kyushu University\\  
6-10-1, Hakozaki, Fukuoka, 812-8581\\ JAPAN}\\ 
\textit{E-mail address}: wakayama@math.kyushu-u.ac.jp\\ 
\end{flushleft}


\begin{thebibliography}{VW}

\bibitem[Ar]{Ar} E. Artin, 
\textit{\"{U}ber die Zetafunktionen gewisser algebraischer Zahlk\"{o}rper}, 
Math. Ann. {\bf 89}(1923), 147--156.

\bibitem[Di]{Di} L. E. Dickson, 
\textit{Linear groups: With an exposition of the Galois field theory},
 Dover Phoenix Editions, Dover Publications, Inc., New York, 1958.

\bibitem[GW]{GW} R. Gangolli and G. Warner, 
\textit{Zeta functions of Selberg's type for some noncompact quotients of symmetric spaces 
of rank one},
Nagoya Math. J. {\bf 78}(1980), 1--44.

\bibitem[G]{G} C. F. Gauss, 
\textit{Disquisitiones arithmeticae}, Fleischer, Leipzig, 1801.

\bibitem[H]{H} Y. Hashimoto,
\textit{Arithmetic expressions of Selberg's zeta functions 
for congruence subgroups}, to appear in J. Numb. Theory.

\bibitem[He]{He} D. Hejhal,
\textit{The Selberg trace formula of $\mathrm{PSL}(2,\Bbb R)$} I, 
Lecture Notes in Math. {\bf 548}, Springer-Verlag, Berlin, 1976/ 
II, Lecture Notes in Math. {\bf 1001}, Springer-Verlag, Berlin, 1983.  

\bibitem[Kl]{Kl} H. D. Kloosterman, 
\textit{The behavior of general theta functions under the modular group 
and the characters of binary modular congruence group, I},
Annals of Math. {\bf 47}(1946), 317--375.

\bibitem[Na]{Na} W. Narkiewicz, 
\textit{Elementary and analytic theory of algebraic numbers}, second edition, Springer-Verlag,
Berlin; PWN, Warsaw, 1990.

\bibitem[Sa]{Sa} P.~Sarnak,
\textit{Class numbers of indefinite binary quadratic forms}, 
J.~Number Theory {\bf 15}(1982), 229--247.

\bibitem[Se]{Se} A.~Selberg,
\textit{Harmonic analysis}, G\"{o}ttingen Lecture Notes (1954), 
Collected papers of A. Selberg vol. 1, 626--674, Springer-Verlag, Berlin, 1989.

\bibitem[Su1]{Su1} T. Sunada,
\textit{Fundamental groups and Laplacians}, Kinokuniya Shoten, Tokyo, 1988 (in Japanese).

\bibitem[Su2]{Su2} T. Sunada, 
\textit{$L$-functions in geometry and some applications}, 
Curvature and topology of Riemannian manifolds (Katata, 1985), 266--284, 
Lecture Notes in Math. {\bf 1201},
Springer, Berlin, 1986.

\bibitem[Ta]{Ta} T. Takagi,
\textit{Algebraic number theory}, Second edition, Iwanami Shoten, Tokyo, 1971 (in Japanese).

\bibitem[Tc]{Tc} N. Tchebotarev, \textit{Die Bestimmung der Dichtigkeit einer Menge von Primzahlen, 
welch zu einer gegebenen Substitutionsklasse gehoren}, Math. Ann. {\bf 95}(1926), 191--228. 

\bibitem[VZ]{VZ} A. B. Venkov and P. G. Zograf, 
\textit{Analogues of Artin's factorization formulas in the spectral theory 
of automorphic functions associated with induced representations of Fuchsian groups},
Math. USSR Izv. {\bf 21}(1983), 435--443.


\end{thebibliography}
\end{document}